\documentclass[a4paper,10pt]{article}
\usepackage[utf8]{inputenc}

\usepackage{amsmath}
\usepackage{multirow}
\usepackage{amssymb}

\usepackage{enumerate}

\newtheorem{thm}{Theorem} 
\newtheorem{theorem}[thm]{Theorem} 
\newtheorem{prop}[thm]{Proposition}

\newtheorem{lemma}[thm]{Lemma} 
 
\newtheorem{cor}[thm]{Corollary}

\newtheorem{defn}[thm]{Definition}

\newcommand{\beq}{\begin{eqnarray}}
\newcommand{\eeq}{\end{eqnarray}}

\newcommand{\Z}{\mathbb Z}
\newcommand{\N}{\mathbb N}

\newcommand{\be}{\begin{enumerate}}
\newcommand{\ee}{\end{enumerate}}\usepackage{enumerate}

\begin{document}
\title{The Multiplicative Automorphisms of a Finite Nearfield, with an Application}
\author{Tim Boykett and Karin-Therese Howell }

\maketitle

\begin{abstract}
In this paper we look at the automorphisms of the multiplicative group of finite nearfields.
We find partial results for the actual automorphism groups.
We find counting techniques for the size of all finite nearfields. 
We then show that these results can be used in order to count the number of near vector spaces of a given dimension over a given nearfield, up to isomorphism.
\end{abstract}

\section{Introduction}

In this paper we look at the automorphisms of the multiplicative group of finite nearfields.
While these are easily found explicitly for fields and the seven exceptional nearfields, we only obtain explicit structural results for some classes of Dickson nearfields.

In the next section we look at finite dimensional near vector spaces.
We show that the results about the multiplicative group automorphisms can be used in order to count the number of near vector spaces of a given dimension over a given nearfield, up to isomorphism.

In the last section we obtain results allowing us to count the size of the automorphism group of a Dickson nearfield, allowing us to calculate the number of  near vector spaces  explicitly.

\section{Background}

A \emph{nearfield} is a $(2,2)$-algebra $(N,+,*)$ such that $(N,+)$ is a group with identity $0$, 
$(N\setminus \{0\},*)$ is a group and one distributive law applies, in our case the right distributive law,
$(a+b)*c = a*c+b*c\,\mbox{ for all } a,b,c \in N$.
Clearly fields are examples of nearfields.
We will often omit the $*$ symbol and notate multiplication by juxtaposition.
The study of nearfields arose from an interest in the axiomatics of fields.
Dickson showed that proper nearfields that are not fields exist using a construction 
that now bears his name, which we will encounter later. 

Nearfields have found important applications in combinatorics and are used in the study of finite geometries, coordinatising certain types of finite planes.
Sharply 2-transitive group actions can be shown to be affine maps $x \mapsto ax+b$ on nearfields in the finite 
case. It has recently been shown that there are non-nearfield based examples in the infinite case \cite{tent}.

By Zassenhaus \cite{zass}, here are three classes that a finite nearfield N can belong to:
\begin{enumerate}
\item $N$ is a finite field,
\item $N$ is a Dickson nearfield,
\item $N$ is one of 7 exceptional nearfields.
\end{enumerate}

In all cases the additive group of a nearfield is elementary abelian. Other facts about nearfields will be introduced below.

We will use $\varphi(n)$ for the Euler Totient function. 
We will use $\Z_n$ for the cyclic group of order $n$, as well as the integers modulo $n$. 
The context should help remove any ambiguity.

\section{Automorphism groups of finite nearfields}

We are interested in exploring the multiplicative structure of finite nearfields through the automorphisms of the group. 
For nearfields that are fields, this is the automorphism group of a cyclic group.
The rest of this section will be spent determining the  multiplicative automorphism group for those finite nearfields that are not fields.

\subsection{Dickson Nearfields}

A Dickson nearfield is a ``twisting'' of a field. The twisting is defined by a Dickson Pair.
\begin{defn}
 A  pair of numbers $(q,n) \in \N^{2}$ is called a \emph{Dickson pair} if $q$ is some power $p^{l}$ of a prime $p$ and each prime divisor of $n$ divides $q-1$ and $4 \vert n \Rightarrow 4 \vert (q-1)$.
\end{defn}

All Dickson nearfields arise by taking Dickson pairs as described in Theorem 8.31, p. 244 in \cite{Pilz}.

For any Dickson pair $(q,n)$ we will denote the associated Dickson nearfield by $DF(q,n)$, the multiplicative group by $G_{q,n}$.
The group is metacyclic and can be presented as 
\[
 G_{q,n} = \left<a,b \,\vert \, a^m=1,\, b^n=a^t,\, ba = a^qb\right>
\]
where $t = m/(q-1)$ and $m = (q^n -1)/n$. \cite[ p 168 (IV.1.1.d)]{waehling}. 

For ease of calculation, we will generally write elements in the form $a^ib^j$.  We note that 
\[(a^kb)^j = a^{k(q^{j-1}+\ldots+q+1)}b^j\] by the third equation above.

The following result gives us some structure of the group $G_{q,n}$, using the numbering of the book, but different notation.

\begin{thm}[\cite{waehling} Theorem IV.1.5]
\label{waehlingthm}
 Let $(q,n)$ be a Dickson pair, $m = (q^n -1)/n$ and $t = m/(q-1)$.
\begin{enumerate}[(a)] 
\item $Z(G_{q,n}) = \langle a^t\rangle$
\item The Sylow subgroups of $G_{q,n}$ are cyclic or quaternion.
\item $\gcd(n,t) = \gcd(q-1,t) \leq 2$ and the following are equivalent:
\begin{enumerate}[(I)]
\item The 2-Sylow subgroups are cyclic
\item $\gcd(q-1,t)=1$
\item $n$ is odd or $q \not\equiv 3 \pmod 4$
\end{enumerate}
\item $\langle a \rangle$ is a maximal abelian subgroup of $G_{q,n}$
\item If $(q,n)\neq (3,2)$ then $\langle a \rangle$ includes all abelian normal subgroups of $G_{q,n}$. In particular, $\langle a \rangle$ is characteristic in $G_{q,n}$.
\end{enumerate}
\end{thm}

The automorphisms of a Dickson nearfield are determined in the following main result.
The multiplicative group $G_{3,2}$ is the quaternion group, which has 24 automorphisms, isomorphic to $S_4$.

\begin{thm}
\label{thm_aut}
 Let $(q,n)$ be a Dickson pair, $(q,n) \neq (3,2)$, with $m=\frac{q^n-1}{n}$ and $t=\frac{m}{q-1}$.
Then  $\phi \in Aut(G_{q,n})$ iff there exist some $i,k \in \Z_m$ such that $\phi(a) = a^i$ and $\phi(b) = a^kb$,
with $\gcd(i,m)=1$ and solving 
\beq
k\frac{q^n-1}{q-1} +t \equiv it \pmod {m} \label{ktnequation}
\eeq
\end{thm}
Proof:
$(\Leftarrow)$
Every map of the generators into $G_{q,n}$ that respects the equations in the presentation can be extended to an endomorphism.
We use $\phi_{ik}$ to denote the morphism with $\phi_{ik}(a) = a^i$ and $\phi_{ik}(b) = a^kb$.
Firstly we show that the images satisfy the three equations of the presentation.
\begin{align}
 \phi_{ik}(a)^m &= a^{im} = 1\\
 \phi_{ik}(b)^n &= (a^kb)^n = a^{k(1+q+\ldots+q^{n-1})} b^n \\
                &= a^{k(1+q+\ldots+q^{n-1})+t} \\
                &= a^{it} = \phi_{ik}(a)^t\\
 \phi_{ik}(b)\phi_{ik}(a) =a^kba^i &= a^ka^{iq}b = a^{iq}a^kb = \phi_{ik}(a)^q\phi_{ik}(b)
\end{align}

Now we show that the homomorphism satisfying these requirements is a bijection.
Suppose $a^\alpha b^\beta \in \ker \phi_{ik}$. 
Then
\begin{align}
 \phi_{ik} (a^\alpha b^\beta) =1 &\Leftrightarrow a^{\alpha i} (a^k b)^\beta = 1 \\ 
&\Leftrightarrow a^{\alpha i} a^{k(1+q+\ldots+q^{\beta-1})} b^\beta = 1.
\end{align}
Thus $\beta=0$, so $\phi_{ik} (a^\alpha) =1 \Leftrightarrow a^{\alpha i} = 1$. 
Because $i$ is coprime to $m$, this implies that $\alpha = 0$ so the kernel of $\phi_{ik}$ is trivial and we have an automorphism.

$(\Rightarrow)$ We show the form of the automorphisms by investigating the images of the generators under the automorphism. 
If mapping of the generators $a,b$ into $G_{q,n}$ is an automorphism 
then the images satisfy the same equations.

By 
Theorem \ref{waehlingthm}
we know that the subgroup generated by $a$ is characteristic, so $a$ must map to some power of itself. 
In order for $\phi$ to be an automorphism, this power must be coprime to the order of $a$, which is $m$, since $a^{m}=1$.

In general, $\phi(b) = a^kb^j$ for some $j,k$, since $ba = a^qb$. The images of $a$ and $b$ must satisfy the same equations as $a$ and $b$, 
so $ba = a^qb$  requires
\begin{align}
 (a^k b^j) a^i &= a^{iq}a^kb^j \\
 \Leftrightarrow a^k a^{i q^j} b^j &= a^ka^{iq}b^j \\
 \Leftrightarrow a^{i q^j} &= a^{iq} \\
 \Leftrightarrow i q^j &\equiv iq \pmod{m}\\
 \Leftrightarrow  q^j &\equiv q \pmod{m}  \mbox{ because $i$ is coprime to $m$}
\end{align}

Let $\alpha = j-1$ and we have the condition $q^\alpha \equiv 1 \pmod{m}$ or equivalently,
$m \vert q^\alpha -1$.
Taking logarithms to the base $q$, we obtain
\begin{align}
 \log_q m &\leq \log_q (q^\alpha -1) \leq \alpha \mbox{ and}\\
 \log_q m &= \log_q \frac{q^n -1}{n} = \log_q (q^n -1) - \log_q n > n - 1 - \log_q n
\end{align}
so $ \alpha > n - 1 - \log_q n$.

The set $A=\{\alpha \in \Z_n \, \vert \, q ^ \alpha \equiv 1 \pmod{m}\} \leq (\Z_n,+)$ as a group, 
and $0 \neq \alpha \in A \Rightarrow n-\alpha \in A$. 
Note also that $A$ is always a proper subgroup, as $1 \not\in A$.

Suppose $A$ is not trivial, so $n \geq 4$. 
Then there is some $\alpha \in A$, $\alpha > n - 1 - \log_q n$.
So $n-\alpha < n-(n-1- \log_q n) = 1+\log_q n$.
But also $n-\alpha \in A$, so $n - \alpha > n - 1 - \log_q n$, so $1+\log_q n > n -1-\log_q n$, that is,
$2 \log_q n > n-2$ so $n^2 > q ^{n-2}$. 
This can be rewritten as $q < n^\frac{2}{n-2}$.
For $n=4$, the upper bound is $ n^\frac{2}{n-2} = 4$ and the bound is monotone descending.
There are no Dickson pairs $(q,n)$ with $q<4$ and $n\geq 4$.
Thus $A$ is always trivial, so we know that $j=1$.

Thus we know that $\phi(b) = a^k b$ for some $k \in \Z_m$.

The second equation in the presentation then requires that 
\begin{align}
 (a^k b)^n &= (a^i)^t \\
 \Leftrightarrow a^{k(q^{n-1}+\ldots +q+1)} b^n &= a^{it} \\
 \Leftrightarrow a^{k(q^{n-1}+\ldots +q+1)} a^t &= a^{it} \\
 \Leftrightarrow k(q^{n-1}+\ldots +q+1)+t &\equiv it \pmod{m}
\end{align}
which is the condition we wanted.
\hfill$\Box$

\medskip

The structure of the automorphism group is a semidirect product.
\begin{lemma}
 Let ${\cal U}(\Z_m)$ denote the group of units of $\Z_m$ under multiplication. Then $Aut(G_{q,n}) \leq {\cal U}(\Z_m) \rtimes \Z_m$, with the units acting by multiplication.
\end{lemma}
Proof:
Let $(i,k),(j,l) \in {\cal U}(\Z_m) \times \Z_m$ and define $(i,k)*(j,l)= (ij,k+il)$.
This gives a multiplication making $({\cal U}(\Z_m) \times \Z_m,*)$ a group that is a semidirect product of the units of $\Z_m$ with $\Z_m$.

Let $\phi_{ik},\phi_{jl} \in Aut(G_{q,n},*)$. Then
\begin{align}
 \phi_{ik} \circ \phi_{jl} (a) &= \phi_{ik}  (a^j) = a^{ij}\\
 \phi_{ik} \circ \phi_{jl} (b) &= \phi_{ik}  (a^lb) = \phi_{ik}  (a^l)\phi_{ik}  (b) = a^{il}a^kb = a^{k+il}b
\end{align}
so $\phi_{ik} \circ \phi_{jl} = \phi_{(ij)(k+il)}$.
Clearly this is the multiplication defined in the semidirect product above, so the automorphisms form a subgroup of this semidirect product.
\hfill$\Box$

\medskip

We note that the inner automorphisms are easily found.

\begin{lemma}
 Let $(q,n)$ be a Dickson pair. Then the inner automorphisms of $G_{q,n}$ are of order $tn$.
\end{lemma}
Proof: The center of the multiplicative group of a Dickson nearfield is of order $q-1$, so 
$\vert Inn(G_{q,n})\vert = \frac{q^n-1}{q-1} = \frac{q^n-1}{n} \frac{n}{q-1}=\frac{m}{q-1} n  = tn$.
\hfill$\Box$

\medskip

We can determine the structure of the automorphism group directly in some cases.

\begin{defn}[\cite{golasinski09}]
 Let $m,n,k \in \N$ such that $k^n\equiv 1 \pmod{m}$. Then $D(m,n;k)$ is the metacyclic group of order $mn$ defined by the presentation
 \begin{align}
  \langle x,y \vert x^m=y^n=1,\, yxy^{-1} = x^k\rangle.
 \end{align}
\end{defn}

For $m,n$ coprime, there exists an explicit description of the automorphism group of $D(m,n;k)$.
\begin{prop}[\cite{golasinski05} Proposition 1.3]
\label{golasinskiprop}
 Let $A$ and $G$ be finite groups of relatively prime orders, $A$ cyclic with a $G$-action $\alpha$.
 Then 
 \[
  Aut(A \times_\alpha G) \equiv A/A^G \times_* (Aut(A) \times Aut_\alpha(G))
 \]
with $Aut_\alpha(G) = \{\phi \in Aut(G) \vert \alpha = \alpha \phi\}$,
$A^G = \{a\in A \vert a^g=a \mbox{ for all } g\in G\}$ and action
$(\phi_1,\phi_2) * x = \phi_1(x)$ for all $(\phi_1,\phi_2) \in Aut(A) \times Aut_\alpha(G)$ and $x \in A/A^G$.
\end{prop}

\begin{thm}
 Let $(q,n)$ be a Dickson pair with $n$ odd or $q \not\equiv 3 \pmod{4}$. 
 Then $G_{q,n} \cong D(r,s;q^{\bar r})$ with $r=t\bar r$, $s=n\bar s$ and $q=\bar r \bar s +1$
\end{thm}
Proof:
The order of $G_{q,n}$ is $q^n-1=mn = t(q-1)n$.
By Theorem \ref{waehlingthm} $\gcd(n,t)=\gcd(q-1,t)=1.$
Let $\bar r \bar s = q-1$ such that all prime factors of $n$ are in $\bar s$, $\gcd(\bar r,n)=1$.
Then define $r=t\bar r$ and $s=n\bar s$ so $rs=q^n-1$.

Let $a,b$ be the generators of $G_{q,n}$.
Define $\bar a = a^{\bar s}$, $\bar b = b^{\bar r}$. Then

\begin{align}
 \bar a^r = a^{\bar s r} = a^{\bar s t \bar r} = a^{t(q-1)}=a^m=1\\
 \bar b^s = b^{\bar r s} = b^{\bar r \bar s n} = a^{t \bar r \bar s} = a^{t(q-1)} = 1\\
 \bar b \bar a \bar b^{-1} = b^{\bar r} a^{\bar s}  \bar b^{-1} = a^{\bar r q^{\bar r}}  \bar b^{-1}
   = \bar a^{q^{\bar r}}\bar b \bar b^{-1} = \bar a^{q^{\bar r}}.
\end{align}

Thus $A =\langle \bar a, \bar b\rangle \leq G_{q,n}$ is a homomorphic image of $D(r,s;q^{\bar r})$.

We now show that $A = G_{q,n}$, thus of the same order as $D$, thus isomorphic.
By calculation, $\bar b^n = a^{t\bar r} = a^r$ and $\gcd(r, \bar s) = 1$ so $a \in \langle a^r,a^{\bar s} \rangle$.
Since $\gcd(\bar r, n) = 1$ there exists some $j$ such that $\bar r j \equiv 1 \pmod n$
so $\bar b ^j = b^{\beta n + 1} = a^{\beta t}b$.
Since $a \in A$ then $a^{\beta t} \in A$ and so $a,b \in A$, $A = G_{q,n}$ and we are done.
\hfill$\Box$

\medskip
Thus we can use Proposition \ref{golasinskiprop} to determine the structure of the automorphism group in this case.
We have not been able to obtain explicit descriptions of the automorphism group, or 
any description of the automorphism group in the case that $\gcd(n,t)=2$. 
However, in Corollary \ref{corrcount} below we will obtain an explicit size of the automorphism group for all finite Dickson nearfields.

\subsection{Exceptional Nearfields}

The nearfield multiplicative groups  are given in \cite[Kapital IV]{waehling}, allowing us to simply calculate the values needed.
The following facts are collected from \cite[(IV.7.1),(IV.8.1)]{waehling} with calculations in \cite{SONATA,GAP4}.
\[
\begin{array}{|c|c|c|c|}
 \hline
\mbox{Name} & \mbox{Mult Gp} & \mbox{Group Aut} \\
\hline
(I) 5^2 & SL(2,3) & S_4 \\  
\hline 
(II) 11^2 & SL(2,3) \times \Z_5 & S_4 \times \Z_4 \\
\hline
(III) 7^2 & 2O & S_4 \times \Z_2 \\ 
\hline
(IV) 23^2 & 2O \times \Z_{11} & S_4 \times \Z_2 \times \Z_{10} \\
\hline
(V) 11^2 & SL(2,5) & S_5 \\ 
\hline
(VI) 29^2 & SL(2,5)\times \Z_7 & S_5 \times \Z_6 \\
\hline
(VII) 59^2 & SL(2,5)\times \Z_{29} & S_5 \times \Z_{28} \\  
\hline
\end{array}
\]

We use $2O$ to represent the binary octahedral group, $SL(n,q)$ the special linear group of dimension $n$ over the finite field of order $q$.

%

\section{An application to counting near vector spaces}

In \cite{Andre}, the concept of a vector space, i.e., linear space, is generalized by Andr\'e to a structure comprising a bit more non-linearity, the so-called near vector space. In \cite{vdWalt2} van der Walt showed how to construct an arbitrary finite-dimensional nearvector space, using a finite number of nearfields, all having isomorphic multiplicative semigroups. In \cite{HowM2}  this construction is used to characterize and count all finite-dimensional nearvector spaces over arbitrary finite fields. 


Our aim is to use our results to count the number of non-isomorphic nearvector spaces over a given nearfield.

\subsection{Near vector spaces}

We begin with a brief overview and some definitions. See \cite{Andre} for further details.

\begin{defn} [\cite{Andre}] A pair $(V,A)$ is called a \emph{near vector space} if:
\be
\item $(V,+)$ is a group and $A$ is a set of endomorphisms of $V$;
\item $A$ contains the endomorphisms $0$, {\it id} and $-${\it id};
\item $A^*=A\setminus\{0\}$ is a subgroup of the group Aut$(V)$;
\item $A$ acts fixed point free (fpf) on $V$, i.e., for $x\in V, \alpha,\beta\in A$, $x\alpha=x\beta$ implies that $x=0$ or $\alpha=\beta$;
\item the quasi-kernel $Q(V)$ of $V$, generates $V$ as a group. Here, $Q(V) = \{x\in V \,|\, \forall \alpha,\beta\in A, \exists\gamma\in A \mbox{ such that } x\alpha + x\beta = x\gamma\}$. \ee
\end{defn}

We sometimes refer to $V$ as a {\it near vector space over $A$}. The elements of $V$ are called {\it vectors} and the members of $A$ {\it scalars}. The action of $A$ on $V$ is called {\it scalar multiplication}. Note that $-${\it id} $\in A$ implies that $(V,+)$ is an abelian group. Also, the dimension of the near vector space, $\dim(V)$, is uniquely determined by the cardinality of an independent generating set for $Q(V)$. 

\medskip

In  \cite{vdWalt2}(Theorem 3.4, p.301) van der Walt derives a characterization of finite dimensional nearvector spaces:

\medskip

\begin{thm} Let $V$ be a group and let $A:=D \,\cup\, \{0\}$, where $D$ is a fix point free group of automorphisms of $V$. Then $(V,A)$ is a finite dimensional near vector space if and only if there exists a finite number of nearfields, $F_{1},F_{2},\ldots,F_{n}$, semigroup isomorphisms $\psi_{i}:A \rightarrow F_{i}$ and a group isomorphism $\Phi:V \rightarrow F_{1}\oplus F_{2}\oplus \cdots \oplus F_{n}$ such that if \[\Phi(v) = (x_{1},x_{2}, \cdots, x_{n}),\,\,\,(x_{i} \in F_{i})\]
then
\[\Phi(v\alpha) = (x_{1}\psi_{1}(\alpha), x_{2}\psi_{2}(\alpha), \cdots, x_{n}\psi_{n}(\alpha)),\]
for all $v \in V$ and $\alpha \in A$.
\end{thm}

\medskip

According to this theorem we can specify a finite dimensional near vector space by taking $n$ nearfields $F_{1},F_{2},\ldots,F_{n}$ for which there are semigroup isomorphisms $\vartheta_{ij}:(F_{j},\,\cdot)\rightarrow (F_{i},\,\cdot)$ with $\vartheta_{ij}\vartheta_{jk}=\vartheta_{ik}$ for $1 \leq i,j,k \leq n$.
We can then take $V := F_{1}\oplus F_{2} \oplus \cdots \oplus F_{n}$ as the additive group of the near vector space and any one of the semigroups ($F_{i_{o}}$, $\cdot$) as the semigroup of endomorphisms by defining
\[(x_{1},x_{2},\ldots,x_{n})\alpha := (x_{1}\vartheta_{1i_{o}}(\alpha), x_{2}\vartheta_{2i_{o}}(\alpha),\cdots, x_{n}\vartheta_{ni_{o}}(\alpha)),\]
for all $x_{j} \in F_{j}$ and all $\alpha \in F_{i_{o}}$. 

\medskip

\begin{defn}
 Two near vector spaces $V$ and $W$ over the same nearfield $N$ are isomorphic if there is a group isomorphism $\theta: (V,+) \rightarrow (W,+)$ such that
$\theta(vn) = \theta(v)n$ for all $v\in V$ and $n\in N$.
\end{defn}

Finite near vector spaces have a finite dimension and a finite base nearfield.
Thus we can represent a finite near vector space over $N$ as a finite sequence of multiplicative automorphisms of $N$, at least one of which is the identity.
An isomorphism between two such representations $(\alpha_1,\ldots,\alpha_n)$ and $(\beta_1,\ldots,\beta_n)$ is a permutation $\gamma$ of $\{1,\ldots,n\}$ and a collection of nearfield automorphisms $\delta_i$ such that 
\[\alpha_i = \delta_i \beta_{\gamma i}\]

\begin{lemma}
 An isomorphism class of near vector spaces of dimension $n$ over a nearfield $N$ is
a multiset of size $n$ of cosets of $Aut(N,+,\cdot)$ in $Aut(N,\cdot)$, at least one of which is $Aut(N,+,\cdot)$.
\end{lemma}
Proof: A multiset is a sequence modulo permutations. So the $\gamma$ above is identity. Then the condition above is
\[ \alpha_i \beta_i^{-1} \in Aut(N,+,\cdot)
\]
\hfill$\Box$

Let us order the cosets as $K_1,K_2,\ldots, K_k$ with $K_1=Aut(N,+,\cdot)$.
Then an isomorphism class representative can be written as a nondecreasing sequence of length $n$, starting with 1, of elements in $1,\ldots,k$.
This number is $k$ multichoose $n-1$, or
\[
 \binom{n+k-2}{n}
\]
by the counting result for multisets \cite[ Theorem 5.3.2]{tucker}.

The rest of this paper will be concerned with calculating the factor $k$ in the above result.

\section{Counting}

This section will be concerned with calculating the size of the automorphism group of a nearfield.
We will look at three cases corresponding to Zassenhaus' result.
We will call the factor  $F(N)$, the index of the nearfield automorphism group within the multiplicative automorphism group of a nearfield $N$.

\subsection{Fields}

Suppose $N=GF(q)$ is a finite field of order $q=p^n$.
Then the field automorphisms are generated by the Frobenius automorphism $x\mapsto x^p$ and so the size of the automorphism group is $n$. 
The multiplicative group of a finite field is cyclic of order $q-1$, so there are $\varphi(q-1)$  generators. 
Thus there are $\varphi(q-1)$ multiplicative automorphisms and the index of the nearfield automorphisms in the multiplicative automorphisms is
\[
 F(GF(q)) = \frac{\varphi(q-1)}{  n} 
\]
This agrees with the result in \cite[Theorem 3.2]{HowM1} for prime fields and \cite[Theorems 3.8 and 3.9]{HowM2} in general.

\subsection{Exceptional Nearfields}

The nearfield multiplicative groups and their automorphism groups 
are given in \cite[Kapital IV]{waehling}, allowing us to simply calculate the values needed.
The following facts are collected from \cite[(IV.7.1),(IV.8.1)]{waehling} with calculations in \cite{SONATA,GAP4}.
\[
\begin{array}{|c|c|c|c|}
 \hline
\mbox{Name} & \mbox{Nearfield Aut} & \mbox{Group Aut} & \mbox{Factor}\\
\hline
(I) 5^2& 4 & 24 & 6\\
\hline 
(II) 11^2 & 2 & 96 & 48\\
\hline
(III) 7^2 & 3 & 48 & 16\\
\hline
(IV) 23^2 & 1 & 480 & 480\\
\hline
(V) 11^2& 5 & 120 & 24\\
\hline
(VI) 29^2 & 2 & 720 & 360\\
\hline
(VII) 59^2 & 1 & 3360 & 3360\\
\hline
\end{array}
\]

So we are done here, and turn our attention to the case of Dickson nearfields.

\subsection{Counting for Dickson Nearfields}

Let $(q,n)$ be a Dickson pair, with $q=p^l$.
The nearfield automorphism group has order 6 for (3,2), otherwise it has order $ln/k$ where $k$ is the order of $p \pmod n$ \cite[Theorem (2.3), p.\ 175]{waehling}.

Our problem now is to count the number of solutions to the equations in Theorem \ref{thm_aut} for a given Dickson pair. 

\begin{defn}
  Let $(q,n)$ be a Dickson pair.
 Define 
 \[
  S(q,n) = \{(i,k) \in \Z_m \times \Z_m \, \vert \, i,k \mbox{ satisfies the equations in Theorem \ref{thm_aut}} \}
 \]
and
\[ T(q,n) = \{(i,k) \in \Z_m \times \Z_{q-1} \, \vert \, \gcd(i,m)=1,\, kn \equiv i-1 \pmod {q-1}\}\]
with $m=\frac{q^n-1}{n}$.
\end{defn}
%

\begin{lemma}
  Let $(q,n)$ be a Dickson pair.
Then $\vert S(q,n)\vert=t \vert T(q,n)\vert$ with $t=\frac{m}{q-1}$.
\end{lemma}
Proof:
The equation (\ref{ktnequation}) above can be simplified. Noting that $t \vert m$ with $m=t(q-1)$ we see that
$k$ satisfies (\ref{ktnequation}) iff $k$ satisfies
\beq
kn \equiv i-1 \pmod {q-1}. \label{knequation}
\eeq
Moreover, if $k$ is a solution for (\ref{ktnequation}) then $k+(q-1)$ is as well. 
Thus one solution to (\ref{knequation}) gives $t$ solutions to (\ref{ktnequation}).
\hfill$\Box$

\begin{lemma}
 Let $(q,n)$ be a Dickson pair.
 For each $i$, the set $\{k \,\vert \, (i,k) \in T(q,n)\}$ is either empty or of order $\gcd(n,q-1)$.
\end{lemma}
Proof:
Suppose $K=\{k \vert (i,k) \in T(q,n)\}$ is nonempty.
Let 
$n=n^* \gcd(n,q-1)$, $q-1 = q^* \gcd(n,q-1)$ so $\gcd(n^*,q^*)=1$.
Then $k \in K$ implies that 
$(k + \frac{q-1}{\gcd(n,q-1)})n = kn + \frac{q-1}{\gcd(n,q-1)}n = kn+(q-1) n^* \equiv kn \pmod{q-1}$ 
so $k + \frac{q-1}{\gcd(n,q-1)}\in K$.
Thus there is an element of $K$ in $\{0,\ldots,\frac{q-1}{\gcd(n,q-1)}-1\}$. 

Suppose there are two. 
Let $k_1,k_2 \in \{0,\ldots,\frac{q-1}{\gcd(n,q-1)}-1\}\cap K$, $k_1 >k_2$. Then
\begin{align}
(k_1 - k_2)n &\equiv (i-1)-(i-1) = 0 \pmod {q-1} \\
 \Rightarrow (q-1) &\vert (k_1-k_2)n \\
 \Rightarrow (k_1-k_2)n &= c(q-1) \,\,\, \mbox{ for some }c \geq 1 \\
 \Rightarrow (k_1-k_2)n^* &= cq^* \\
  k_1,k_2 < q^* &\mbox{ so } (k_1-k_2) < q^* \Rightarrow n^* > c 
\end{align}
But $n^* \vert c$ by coprimeness of $n^*$ and $q^*$, which is a contradiction. 
Thus we cannot have two distinct $k_1,k_2<q^*$ in $K$. 
Thus if $K$ is nonempty, it has precisely $\gcd(n,q-1)$ elements. 
Writing $K$ as a subset of $\Z_{q-1}$, $K$ is a coset of $\langle q^*\rangle \leq \Z_{q-1}$.
\hfill$\Box$

\medskip Thus we can calculate the order of our multiplicative automorphism group.

\begin{cor}\label{corrho}
 Let $(q,n)$ be a Dickson pair.
 Then 
 \begin{align}
 \vert Aut(DF(q,n),*)\vert = \frac{1}{\rho}\,\varphi(m)t \gcd(n,q-1)
 \end{align}
 where $\frac{1}{\rho}$ is the proportion of
 $i$ coprime to $m$ for which $\{k \vert (i,k) \in T(q,n)\}$ is nonempty.
\end{cor}
Proof:
Let $A= \{i \vert (i,k) \in S(q,n)\}$. This is a subgroup of ${\cal U}(\Z_m)$ of index $\rho$. Then by the above two results, 
\begin{align}
 \vert Aut(G_{q,n})\vert &= \vert A \vert t \,\gcd(n,q-1) \\
    &= \frac{1}{\rho} \varphi(m) t\,\gcd(n,q-1).
\end{align}
\hfill$\Box$

We can find an explicit description of the factor $\rho$.

\begin{theorem}
\label{thmrho}
 Let $(q,n)$ be a Dickson pair, $q_2$ the factor of $q-1$ containing all the primes in $n$.
Then $\rho = \frac{\varphi(q_2)\gcd(n,q-1)}{q_2}.$
\end{theorem}
Proof:
First we show that there exists a  $k$ such that $(i,k)\in S(q,n)$ iff $\gcd(i,m)=1$ and $\gcd(n,q-1)\vert (i-1)$.
We know that there exists a $k$ such that $(i,k)\in S(q,n)$ iff $\gcd(i,m)=1$ and $kn\equiv i-1 \pmod{(q-1)}$, so we only need to show
the equivalence of the second requirements.
Suppose  $kn\equiv i-1 \pmod{(q-1)}$, i.e.\ $kn-(i-1)\equiv 0\pmod{(q-1)}$.
Let $g=\gcd(n,q-1)$, $n=\bar n g$ and $q-1=\bar q g$.
Then
\begin{align}
 q-1 \vert kn-(i-1) &\Leftrightarrow \bar q g \vert k\bar n g - (i-1) \\
 & \Rightarrow g \vert (i-1).
\end{align}
The other direction requires us to show that divisibility suffices.
Let $g \vert i-1$, so $\bar i g=i-1$ for some integer $\bar i$.
Then
$k\bar n g - \bar i g \equiv 0 \pmod {\bar q g} \Leftrightarrow k\bar n - \bar i \equiv 0 \pmod {\bar q}$.
We know that $\bar n$ is coprime to $\bar q$ so $\bar n$ is a unit in $\Z_{\bar q}$.
Thus this equation is solvable with $k = \frac{\bar i}{\bar n}$ and $(i,k) \in S(q,n)$, so we have the first part of the proof.

Let $A=\{i\in \Z_m\vert \mbox{there exists a } k \mbox{ such that }(i,k) \in S(q,n)$ with $S(q,n)= 
\{i\in \Z_m\vert \gcd(i,m)=1,\,g\vert (i-1)\}$.
We want to know the size of $A$.

There are two cases to consider, depending upon the value of $\gcd(n,t)$.

Case $\gcd(n,t)=\gcd(q-1,t) = 1$. 
Let $q-1=q_1q_2$ where all the primes in $n$ are in $q_2$, so $\gcd(q_1,n)=1$.
We know $m=t(q-1) = (tq_1)q_2$ with $tq_1$ coprime to $q_2$.
Thus we can write $A=\{i \in \Z_m \vert  \gcd(i,tq_1)=1,\, \gcd(i,q_2)=1,\,g \vert i-1\}$.
For the third condition, we can write $i=\alpha g+1$ for some $\alpha$. 
Because the prime factors of $n$ and $q_2$ are identical, $i$ is then coprime to $q_2$.
Thus we can write
\begin{align}
 A & =\{\alpha g+1 \vert  \gcd(\alpha g+1,tq_1)=1,\, \alpha \in \{0,\ldots,\frac{m}{g}-1\}\}\\
  \bar A & =\{\alpha g+1 \vert  \alpha \in \{0,\ldots,\frac{m}{g}-1\}\}\\
 &= q_2 \Z_m + \{0,\ldots,\frac{q_2}{g}-1\}g + 1
\end{align}
when we write sums of sets as the set of all sums.
For each $\beta \in \{0,\ldots,\frac{q_2}{g}-1\}$, let $A_\beta =  q_2 \Z_m + \beta g + 1$. 
Note that all $A_\beta$ are pairwise distinct and that $\vert A_\beta \vert =\vert q_2 \Z_m\vert = tq_1$.
Then
\begin{align}
 \vert \{i \in A_\beta \vert \gcd(i,tq_1)=1\}\vert = \varphi(tq_1).
\end{align}
Thus $\vert A\vert = \frac{q_2}{ g}\varphi(tq_1)$ so 
\begin{align}
 \rho &= \frac{\varphi(m)}{\vert A\vert}\\ 
 & = \frac{\varphi(tq_1)\varphi(q_2) g}{q_2\varphi(tq_1)}\\ 
 & = \frac{\varphi(q_2) g}{q_2}.
\end{align}

Case $\gcd(n,t) = 2$. The proof is analagous to the proof above, but with changes for the extra factors of 2.
By Theorem \ref{waehlingthm} we know that $n$ is even and $q  \equiv 3 \pmod 4$.
Thus $4$ does not divide $q-1$, so $4$ does not divide $n$.
Let $q-1=q_1q_2$ where all the primes in $n$ are in $q_2$, so $\gcd(q_1,n)=1$.
We know $m=t(q-1) = (\bar tq_1)(2^\tau q_2)$ with $\bar tq_1$ coprime to $2^\tau q_2$, $t=2^\tau \bar t$.
Thus we can write $A=\{i \in \Z_m \vert  \gcd(i,tq_1)=1,\, \gcd(i,2^\tau q_2)=1,\, g \vert i-1\}$.
For the third condition, we can write $i=\alpha  g+1$. 
Because the prime factors of $n$ and $q_2$ are identical, $i$ is then coprime to $2^\tau q_2$.
Thus we can write
\begin{align}
 A & =\{\alpha  g+1 \vert  \gcd(\alpha g+1,\bar t q_1)=1,\, \alpha \in \{0,\ldots,\frac{m}{ g}-1\}\}\\
  \bar A & =\{\alpha  g+1 \vert  \alpha \in \{0,\ldots,\frac{m}{ g}-1\}\}\\
 &= (2^\tau q_2) \Z_m + \{0,\ldots,\frac{2^\tau q_2}{ g}-1\} g + 1
\end{align}
For each $\beta \in \{0,\ldots,\frac{2^\tau q_2}{ g}-1\}$, let $A_\beta =  2^\tau q_2 \Z_m + \beta g + 1$. 
Note that all $A_\beta$ are pairwise distinct and that $\vert A_\beta \vert =\vert 2^\tau q_2 \Z_m\vert = \bar tq_1$.
Then
\begin{align}
 \vert \{i \in A_\beta \vert \gcd(i,\bar tq_1)=1\}\vert = \varphi(\bar tq_1).
\end{align}
Thus $\vert A\vert = \frac{2^\tau q_2}{ g}\varphi(\bar tq_1)$ so 
\begin{align}
 \rho &= \frac{\varphi(m)}{\vert A\vert}\\ 
 & = \frac{\varphi(\bar tq_1)\varphi(2^\tau q_2) g}{2^\tau q_2\varphi(\bar tq_1)}\\ 
 & = \frac{\varphi(2^\tau q_2) g}{2^\tau q_2}\\
 & = \frac{\varphi(2^{\tau+1})\varphi(q_2/2) g}{2^\tau q_2}\\
 & = \frac{2^{\tau}\varphi(q_2/2) g}{2^\tau q_2}\\
 & = \frac{\varphi(q_2/2) g}{q_2}\\
 & = \frac{\varphi(q_2) g}{q_2}
\end{align}
because $s$ odd implies that $\varphi(s) = \varphi(2s)$.

Thus we are done.
\hfill$\Box$

By combining the expressions above, we obtain the following explicit counting result.

\begin{cor}\label{corrcount}
 The size of the automorphism group of $G_{q,n}$ is
 \begin{align}
  t \varphi(t) \varphi(q_1) q_2
 \end{align}
where $t = \frac{q^n-1}{n(q-1)}$ and $q-1 = q_1q_2$ with all primes factors of $n$ in $q_2$.
\end{cor}

\subsection{Applying counting structure}

While the above results give us an explicit expression for a given Dickson pair, this is hard to calculate in general. 
However;  we can show that for certain types of Dickson pairs, there are simpler expressions for the size of the automorphism group and the index of the nearfield automorphism group in it.

\begin{lemma}
 Let $(q,2)$ be a Dickson pair, $q=p^l$. 
 Then $\vert S(q,2)\vert = \varphi(\frac{q^2-1}{2})(q+1)$, so 
 $F(DF(q,2)) = \frac{(q+1)}{2l}\varphi(\frac{q^2-1}{2}).$
\end{lemma}

Proof:  Following the notation above, $q_2 = 2^\tau$ for some $\tau$.
By Theorem 17 we have that 
$\rho = \frac{\varphi(q_2)\gcd(n,q-1)}{q_2} = \frac{2^{\tau-1}}{2^\tau}\gcd(2,q-1)= \frac{\gcd(2,q-1)}{2}$. 
Now using Corollary 16 we get that
\begin{align}
Aut(DF(q,n),*)\vert &= \frac{1}{\rho}\,\varphi(m)t \gcd(2,q-1)\\
 & = \frac{2}{\gcd(2,q-1)}\,\varphi\big(\frac{q^{2}-1}{2} \big)t \gcd(2,q-1)\\\\ 
 & = \varphi\big(\frac{q^{2}-1}{2} \big)(q+1)
\end{align}
As we saw above, the order of the nearfield automorphism group is $ln/k$, where $k$ is the order of $p$ modulo $n$. 
But $n=2$ and $p$ is odd, 
so $k=1$. So the nearfield automorphism group has order $2l$. Thus our factor is the order of $S(q,2)$ divided by $2l$.
\hfill$\Box$

Similar results can be obtained for $n=3$ and others.

\section{Summary}

We have investigated the structure of the multiplicative group of a finite nearfield and found explicit counts for the Dickson nearfield case, as well as the seven exceptional finite nearfields.

These results can be used to count the number of nonisomorphic near vector spaces of a given dimension over a given finite nearfield.
We showed that finite near vector spaces have isomorphism classes that can be readily defined by sequences of cosets. 
Thus the process of enumerating examples is based upon questions of the multiplicative and nearfield automorphism groups and the index of the latter in the former.

The main open problems remaining here relate to the explicit construction of the automorphism groups. 
Numerically we were able to find the size of the automorphism group in all cases, using some
number theoretical properties.
We were able to show that in the case that $\gcd(n,t)=1$, the multiplicative group of a Dickson nearfield
is a split metacyclic group and thus we can in principle determine the automorphism group. However we were
unable to find an explicit expression for the automorphism group. 
It would appear that techniques similar to those used in Theorem \ref{thmrho} could be applied to determine the automorphism group explicitly.

Some smaller questions require attention. 
Can we explicitly describe the inner automorphisms? Does this tell us something important about the automorphism group?

These results show that there are a large number of nonisomorphic near vector spaces of a given dimension over a given nearfield. 
Some applications of near vector spaces have been proposed \cite{Andre} and this quantity of examples offers some indication of their rich structure.
In \cite{boykettwendt15} one author and Gerhard Wendt have shown that near vector space type constructions, in particular homeomorphisms of near vector spaces, can be used to construct nearrings with special properties, generalising the properties of units acting as endomorphisms on the additive group of a nearfield. 
This is an enticing direction of further research.

\section*{Acknowledgements}

We would like to thank Peter Mayr and the referee for some very useful suggestions. This work has been partially supported by Austrian Science Funding (FWF) Project number 23689-N18 and by the South African National Research Foundation.

\end{document}